\documentclass [12pt] {article}
\usepackage {amsmath}
\usepackage [cp866] {inputenc}
\usepackage [english, russian] {babel}
\usepackage {amsfonts, amssymb}
\usepackage {latexsym}

\newtheorem {theorem} {Теорема} [section]
\newtheorem {lemma} {Лемма} [section]
\newtheorem {cor} {Следствие} [section]
\title {О повторяющейся концентрации и периодических режимах при аномальной диффузии в полимерах}
\author {Д. А. Воротников\footnote{E-mail:
mitvorot@math.vsu.ru; Работа поддержана РФФИ и программой BRHE
Минобрнауки РФ и АФГИР}}
\date {}

\begin {document} \large

\maketitle

УДК 517.9

\begin{abstract} Распространение проникающей жидкости в полимерном материале часто не удовлетворяет классическим уравнениям диффузии и требует учета релаксационных (вязкоупругих) свойств полимера. Мы изучаем краевую задачу на ограниченной области в пространстве для системы уравнений, моделирующих такую аномальную диффузию. Доказано, что для любого доста\-точно короткого отрезка времени и заданного напряжения в начале этого отрезка существует такое глобальное по времени слабое решение краевой задачи (пара концентрация -- напряжение), что концентрация в начале и в конце отрезка одинаковая. При дополнительном условии на коэффициенты показано также существование периодических по времени слабых решений (без ограничений на длину периода).
\end{abstract}

\begin{center}
Ключевые слова: нефиковская диффузия, полимер, проникающая жидкость, топологическая степень, слабое решение, периодичность
\end{center}

\section {Введение}

\numberwithin{equation}{section}

\newcommand {\R} {\mathbb{R}}

\def\be{\begin{equation}}
\def\ee{\end{equation}}
\def\fr#1#2{\frac{\partial #1}{\partial #2}}

Благодаря своим разнообразным и многосторонним свойствам полимеры
находят широкое применение во многих отраслях промышленности. В то
же время строгое научное (в том числе математическое) изучение их
свойств даже отстает от темпов внедрения этих материалов, хотя и
ведется достаточно широко.

Полимеры часто используются в комбинации с жидкостью, проникающей\\ внутрь его структуры.
При этом меняются свойства самого полимера.
Например, в фарма\-цевтической промышленности жидкое действующее вещество может приме\-няться в сочетании с полимерным носителем.
Оказывается, что в такой ситуации процесс распространения жидкости в полимере противоречит законам классичес\-кой теории диффузии. По всей видимости, это объясняется \cite{chn3} сложной, комбини\-рованной, вязкоупругой природой полимеров.

Напомним классический общий закон сохранения \be\fr u t = - div\, J, \ee который, в частности, подхо\-дит для описания диффузии в сплошных средах: в этом случае
$u=u(t,x)$ обозначает концентрацию, а $J=J(t,x)$ обозначает поток концентрации (они зависят от времени $t$ и простран\-ственной точки $x$).
\textit{Закон Фика} гласит, что поток концентрации пропорционален ее градиенту: \be J= -D(u) \nabla u, \ee где $D(u)$
есть коэффициент диффузии (в общей постановке это положительно-определенный тензор; в более простой он вырождается в скаляр и/или перестает зависеть от концентрации). Из (1.1) и (1.2) следует классическое уравнение диф\-фузии: \be\fr u t = div\, (D(u) \nabla
u).\ee

При распространении жидкости в полимере возникает ряд явлений (например, называемые в англоязычной литературе \textit{case II diffusion} \cite{tw,tw1}, \textit{sorption overshoot} \cite{chn1}, \textit{literal skinning, trapping skinning} \cite{ed,ed2} и \textit{desorption overshoot} \cite{cc}), которые показывают, что поведение концентрации жидкости не вписывается в ограничения, налагаемые (1.2) или (1.3). Как и ньютоновское реологическое соот\-но\-шение (см. \cite{gruy}), подходящее для описания жидкостей, но слишком простое для описания полимеров, закон Фика также должен быть значительно изменен для описания этих (и других потенциально возможных при распространении жидкости в поли\-мере) феноменов.
Используя аналогии с реологией вязкоупругих сред и результаты исследований и опытов многих специалистов, Коэн и др. \cite{chn3,chn1,chn0} предложили следующий определяющий закон (кстати, учитывающий и конвек\-цию):

\be J= -D(u)\nabla u - E(u) \nabla \sigma+M(u,\sigma) u, \ee
\be \fr
{\sigma} t + \beta(u) \sigma = f(u, \fr u t). \ee

Здесь $\sigma$ --- вспомогательная переменная, связанная с напряжением  \cite{sw}, вслед\-ствие чего ее саму называют этим термином (англ. \textit{stress}). Известными (опреде\-ляемыми свойст\-вами материалов и внешними условиями) считаются коэффици\-енты: скаляр $\beta$ (по физическому смыслу обратно пропорциональный времени релаксации и зависящий от концентрации), положи\-тельно-определенные тензоры $D$ и $E$ (зави\-сящие от концентрации), скаляр $f$ и отвечающий за конвекцию вектор $M$ (завися\-щие от двух аргументов).  Вот типичный вид функции $\beta$  \cite{chn1}: $$\beta(u)= \frac 1 2 (\beta_R+
\beta_G)+ \frac 1 2 (\beta_R- \beta_G) \tanh (\frac {u- u_{RG}}
{\delta}),$$ где  $\beta_R,\beta_G, \delta, u_{RG}$ --- положительные константы, $\beta_R>\beta_G$;
$D$ обычно имеет подобную (т.е. монотонную с областью резкого возрастания) зависимость от $u$ \cite{cc}.

Согласно  \cite{chn0,sw}, хороший пример функции $E$ --- линейная зависимость от $u$. Впрочем, $E$ часто считают константой --- это считается допустимым упрощением  модели \cite{chn0}.
В то же время еще в \cite{cox} (численно) было показано, что $E$ не должно быть ненулевой константой ($E\equiv 0$ это фактически закон Фика), так как при $E(0)\neq 0$ концентрация $u$ может становиться отрицательной (что не имеет физиче\-ского смысла). С другой стороны, строго доказано (\cite{am2}, см. также \cite{bei}), что при $E(0)= 0$ концентрация $u$ остается неотрицательной, если она была таковой в каждой точке в начальный момент времени и граничные условия неотрицательны. В настоящей статье (конец пункта 3) приводятся соображения в пользу целесооб\-разности условия $E(1)= 0$, и тогда типичным примером $E$ будет функция вида \be E(u)=\frac {\alpha_1 u(u-1)^2}{\alpha_2 +(u-1)^2}\ee c малым $\alpha_2$, а $\alpha_1$ порядка $1$.

Другие примеры $\beta,D$ и $E$ см. в \cite{chn0,cc}. А для функции $f$ типичен следующий вид \cite{chn0}:\be
f(u,u')=\mu(u) u+ \nu(u) u' ,\ee где функции $\mu$ и $\nu$ принимают значения в ограниченном отрезке $[0,a]$.

Из (1.4) и (1.1) следует уравнение аномальной диффузии: \be\fr u t = div(D \nabla u + E\nabla \sigma -Mu).\ee

С более общей точки зрения \cite{am2}, векторный коэффициент $M$, скаляр $\beta$ и тензоры $D$, $E$ могут зависеть от $t,x,u$ и $\sigma$ (см. еще \cite{ed2}).

Начально-краевые задачи для системы (1.8), (1.5) изучались несколькими авто\-рами. В \cite{am2}  показано существование максимального (не глобального по времени) решения.
Глобальная разрешимость задачи представлена в
\cite{am1} при $f=\mu u$, $M\equiv 0$ и не зависящих от концентрации скалярах $D=E$. В этой работе рассмотрен одномерный случай
($0<x<1$), но по всей видимости он обобщается на случай
$x\in \Omega\subset \mathbb{R}^n$ (ограниченная область с гладкой границей).
В
\cite{bei} же доказано  глобальное существование решений при непостоянном $E$ (удовлетворяющем указан\-ному выше условию $E(0)= 0$ и ряду других условий, исключающих, кстати, возможность того, что $E(1)= 0$). При этом предполагается ограниченность снизу положительной константой начальных и граничных условий на концентрацию, что влечет строгую положительность решений, и не дает возможности рассмотреть наличие "сухих" участ\-ков. Большая часть статьи \cite{bei} рассматривает одномерную ситуацию $0<x<1$, но в конце приведена схема переноса результатов на много\-мерные области. Наличие слабых решений начально-краевой задачи в многомерной ограни\-ченной области изучалось в \cite{p1,p2}.
Для неограниченной области $\Omega=\mathbb{R}^n$ известно глобальное существование диссипативных (по сути очень слабых) решений с посто\-янными скалярными $D$ и $E$, и с $M\equiv 0$ (\cite{diss}; использованная там техника примени\-ма и для
$\Omega\subset \mathbb{R}^n$).
Несколько отличная от (1.8), (1.5) модель изучалась в \cite{riv}.

Настоящая статья посвящена проблеме нахождения периодических по времени решений краевой задачи для системы
(1.8), (1.5) и сходной с ней задачей поиска решений с повторяющейся концентрацией при заданном начальном напряжении. Идея поиска решений эволюционных уравнений в частных производных, повторяю\-щихся в два заданных момента времени (то есть обладающих так называемым \textit{reproductive property}, в формульной записи $y(t_1)=y(t_2)$), идет от работы \cite{kan}. В отличие от задачи о периодических решениях, при этом не требуется периодичность коэффициентов и правых частей уравнений. У задачи с одной неизвестной функ\-цией как раз при периодичности коэффициентов и правых частей решение с \textit{reproductive property} порождает периодическое решение.  Мы же в пункте 3 этой работы ищем решение, обладающее \textit{reproductive property} по концентрации и начальному условию по нап\-ряжению (что однако не влечет автоматически наличия периоди\-ческих решений). В пункте 4 мы докажем основной результат пункта 3 (теорему 3.1) о существовании слабого решения этой задачи. В пункте 5 же мы уже иссле\-дуем задачу о периодических решениях, но слабая разрешимость (теорема 5.1) получается при дополнительных условиях на коэффициенты.

\section {Обозначения}
Мы используем стандартную запись $L_p (\Omega) $, $W_p^{m} (\Omega)
$, $H^{m} (\Omega) =$ $W_2^{m} (\Omega) $ $(m\in\mathbb{Z}, 1 \leq
p \leq \infty)$, $H^{m}_0 (\Omega) = \stackrel{\circ}{W}{}_2^{m} (\Omega) $ $(m \in \mathbb{N})$ для пространств Лебега и Соболева на открытом множестве (области)  $\Omega\subset
\R^n$, $n\in \mathbb{N}$.

Скалярное произведение и евклидова норма в
$L_2(\Omega)^k=L_2(\Omega, \R^k)$ обозначаются $(u,v)$ и $\|u\|$
соответственно ($k=$  $1$ или $n$). В пространстве $H^1_0(\Omega)$
используются следующие скалярное произведение и норма:
$(u,v)_1=(\nabla u, \nabla v), \|u\|_1=\|\nabla u\|$. Напоминаем
неравенство Фридрихса \be \|u\|\leq K_\Omega\|u\|_1.\ee
Аналогично в $H^2_0(\Omega)$ используются скалярное произведение и норма\\
$(u,v)_2=(\Delta u, \Delta v), \|u\|_2=\|\Delta u\|$.

Как обычно, пространство $H^{-m}(\Omega)$, $m=1,2$ отождествляется с сопряженным к
$H_0^m(\Omega)$. Действие функционала из
$H ^ {-m}(\Omega) $ на элемент из $H_0^m(\Omega) $ обозначается
$ \langle\cdot, \cdot\rangle $. Напомним, что $\|\varphi\|_{-m}=\sup_{\|w\|_m=1} |\langle\varphi,
w\rangle|.$

Иногда будем писать просто $L_p $, $H^m$ вместо $L_p (\Omega)^k, H^m (\Omega)^k$
и др., $k=1,n$.

Оператор Лапласа $\Delta: H^1_0(\Omega)\to H^{-1}(\Omega)$ является изоморфизмом.
Поэтому оператор \be\Delta^{-1}: H^{-1}(\Omega)\to H^{1}_0(\Omega)\ee тоже является изоморфизмом.
Обозначим $X=X(\Omega)=\Delta^{-1}(H^{1}_0(\Omega))$. Скалярное произведение и норму в $X$ зададим по формулам $(u,v)_X = (\Delta u, \Delta v)_1$, $\|u\|_X = \|\Delta u\|_1$. Обратим внимание, что между $H^{-1}(\Omega)$ и $X(\Omega)$ имеется следующая двой\-ственность:
\be\left\langle u,v\right\rangle _1=-\left\langle u,\Delta v\right\rangle,\ u\in H^{-1}, v\in X,\ee
причем $\left\langle u,v\right\rangle _1=(u,v)_1$ для $ u\in H^1_0, v\in X$.

Записи типа $C (\mathcal{J}; E) $, $C_w (\mathcal{J}; E) $, $L_2
(\mathcal{J}; E) $, $L_{2,loc}
(\mathcal{J}; E) $ и.т.п. используются для обозначения пространств функций на интервале
$\mathcal{J}\subset \mathbb {R} $ со значениями в банаховом пространстве
$E $ (в данном случае соответственно непрерывных, слабо непрерывных, суммируемых с квадратом, локально суммируемых с квадратом). Напомним, что функция $u:
\mathcal{J} \rightarrow E$ называется \textit{слабо непрерывной}, если  для любого непрерывного линейного функционала $g$ на $E$ функция $g(
u(\cdot)): \mathcal{J}\to \mathbb{R}$ непрерывна.

В случае, когда $E$ является функциональным пространством (напр. $L_2(\Omega),\\ H^m(\Omega)$), мы отождествляем элементы $C (\mathcal{J}; E) $,
$L_2(\mathcal{J};E)$ и др. с числовыми функ\-циями, определенными на
$\mathcal{J}\times \Omega$, по формуле
$$u(t)(x)=u(t,x),\, t\in \mathcal{J}, x\in \Omega.$$

Ниже используются также следующие пространства ($T>0$ --- число):
$$ W= W(\Omega, T) = \{\tau\in L_2 (0, T; H^1_0 (\Omega)), \
\tau' \in L_2 (0, T;
 H ^ {-1} (\Omega)) \}$$ $$\|\tau\|_{W}=\|\tau\|_{L_2 (0, T; H^1_0(\Omega))}+\|\tau'\|_{L_2 (0, T; H ^ {-1}(\Omega))};$$

$$ W_1= W_1(\Omega, T) = \{\tau\in L_2 (0, T; X (\Omega)), \
\tau' \in L_2 (0, T;
 H ^ {-1} (\Omega)) \}$$ $$\|\tau\|_{W_1}=\|\tau\|_{L_2 (0, T; X(\Omega))}+\|\tau'\|_{L_2 (0, T; H ^ {-1}(\Omega))};$$

 $$ W_2= W_2(\Omega, T) = \{\tau\in L_2 (0, T; H^2_0 (\Omega)), \
\tau' \in L_2 (0, T;
 H ^ {-2} (\Omega)) \}$$ $$\|\tau\|_{W_2}=\|\tau\|_{L_2 (0, T; H^2_0(\Omega))}+\|\tau'\|_{L_2 (0, T; H ^ {-2}(\Omega))};$$
\cite[лемма III.1.2]{tem} влечет непрерывность вложений $W,W_2\subset
C([0,T];L_2(\Omega))$, $W_1\subset
C([0,T];H^1_0(\Omega))$ (см. также \cite{ggz,gruy}).

Обозначение $|\cdot|$ используется для модуля числа, евклидовой нормы в $\R^n$ и в следующем случае. Пусть
$\R^{n \times n}$ --- пространство матриц порядка $n\times n$ с нормой
$$|A|=\max_{\xi\in \R^n,|\xi|=1} |A\xi|.$$
Символ $\R^{n \times n}_+ \subset$ $\R^{n \times n}$ это множество таких матриц $A$, что
$$ (A\xi, \xi)_{\R^n} \geq d(A)(\xi, \xi)_{\R^n}$$ для некоторого $d(A)\geq 0$ и всех $\xi\in \R^n$.

\section {Задача о повторяемости концентрации}
Пусть полимер заполняет область $\Omega\subset\mathbb{R}^n$,
$n\in \mathbb{N}$. Наиболее важны частные случаи $n=2$
(диффузия в пленках) и $n=3$. Рассмотрим краевую задачу для системы уравнений, описывающей распространение проникающей жидкости в этом полимере:

\be \fr u t = div[D_0(t,x,u,\sigma) \nabla u $$ $$+ E_0(t,x,u,\sigma) \nabla \sigma -M_0(t,x,u,\sigma) u],\ (t,x)\in [0,T]\times \Omega,\ee \be \fr {\sigma} t
+ \beta_0 (t,x,u,\sigma) \sigma = \mu_0(u) u + \nu_0(u) \fr u t,\ (t,x)\in [0,T]\times \Omega, \ee \be u(t,x)=\varphi(t,x),\ (t,x)\in [0,T]\times \partial\Omega. \ee

Здесь $u=u(t,x):[0,\infty)\times \overline{\Omega} \to \R$ это неизвестная концентрация жидкости (в точке $x$ в момент времени $t$), $\sigma=\sigma(t,x):[0,\infty)\times \overline{\Omega}
\to \R$ --- неизвестное напряжение, $\varphi:[0,\infty)\times \partial\Omega \to \R$ --- заданное граничное условие, $\mu_0,\nu_0:\R\to \R$, $D_0$, $E_0: \R^{n+3}=\R\times\R^n\times\R\times\R\to \R^{n\times n}_+,\ \ $ $\beta_0: \R^{n+3}\to \R$, $M_0: \R^{n+3}\to \R^n$ --- известные функции, $\nu_0(\cdot)\geq 0$.

Будем искать решения, имеющие одинаковые распределения концентрации в моменты $0$ и $T>0$, то есть \be u(0,x) =u (T,x),\ x\in\Omega\ee
(без ограничения общности мы считаем начальным нулевой момент времени).

Предположим на время, что область $\Omega$ и участвующие в уравнениях функции достаточно "хорошие"\footnote{Например, функции $\varphi$, $\mu_0,\nu_0, D_0, E_0,\beta_0,M_0$ --- $C^2$-гладкие и ограниченные, граница области также $C^2$-гладкая, и область локально расположена по одну сторону своей границы. Хотя данные преположения достаточно естественны для изучаемой модели (см. Введение, \cite{chn0}, а также конец этого пункта по поводу ограниченности $E_0$ и $M_0$), ниже мы существенно их ослабим.}.  Функцию $\varphi$ будем считать определенной на $[0,\infty)\times\overline{\Omega}$.
При этом, разумеется, необходимо считать, что
\be \varphi | _ {t=0} =\varphi | _ {t=T}.\ee

Введем новую переменную $$\varsigma(\cdot)=\sigma(\cdot)-\int\limits_0^{u(\cdot)}\nu_0(y)\;dy.$$
По сути это "чисто нефиковская" составляющая напряжения. При $\nu_0\equiv 0$ она совпадает с напряжением.

Введем еще новые функции
$$\gamma(\cdot,\cdot,u,\varsigma)=\mu_0(u)-\frac{\beta_0\left(\cdot,\cdot,u,\varsigma+\int\limits_0^u \nu_0(y)dy\right)\int\limits_0^{u(\cdot)}\nu_0(y)\;dy }{u(\cdot)},$$ $$D_1(t,x,u,\varsigma)=$$ $$D_0\left(t,x,u,\varsigma+\int\limits_0^u \nu_0(y)dy\right)+\nu_0(u)E_0\left(t,x,u,\varsigma+\int\limits_0^u \nu_0(y)dy\right)\in \R^{n\times n}_+ ,$$ $$E_1(\cdot,\cdot,u,\varsigma)=E_0\left(\cdot,\cdot,u,\varsigma+\int\limits_0^u \nu_0(y)dy\right),$$ $$M_1(\cdot,\cdot,u,\varsigma)=-M_0\left(\cdot,\cdot,u,\varsigma+\int\limits_0^u \nu_0(y)dy\right),$$ $$\beta_1(\cdot,\cdot,u,\varsigma)=-\beta_0\left(\cdot,\cdot,u,\varsigma+\int\limits_0^u \nu_0(y)dy\right).$$

Перепишем (3.1) и (3.2) в следующем виде:

\be \fr u t = div[D_1(t,x,u,\varsigma) \nabla u + E_1(t,x,u,\varsigma) \nabla \varsigma + M_1(t,x,u,\varsigma) u], \ee \be \fr {\varsigma} t =
\beta_1 (t,x,u,\varsigma) \varsigma + \gamma(t,x,u,\varsigma) u.\ee

Нам не хватает еще одного "временного" условия. Пусть начальное "чисто\\ нефиковское" напряжение также известно:

\be \varsigma(0,x) = \varsigma_0(x),\ x\in\Omega.\ee

При каждом фиксированном $x\in \overline{\Omega}$ рассмотрим задачу Коши
\be \fr {\psi} t =
\beta_1 (t,x,\varphi(t,x),\psi) \psi + \gamma(t,x,\varphi(t,x),\psi) \varphi(t,x), \ee \be \psi| _
{t=0} = \varsigma_0.\ee

Если $\varsigma_0$, $\varphi$, $\beta_1$ и $\gamma$ ограничены, то ее решение $\psi(t,x)$ также априори ограничено на конечных отрезках (ср. \cite{p1,p2}), и потому существует и единственно на всей положительной полуоси. Заметим, что $\varsigma|_{\partial\Omega}$=$\psi|_{\partial\Omega}$.

Подействуем оператором Лапласа на обе части (3.7):
\be \Delta\varsigma '=div[\nabla
(\beta_1 (t,x,u,\varsigma) \varsigma) + \nabla(\gamma(t,x,u,\varsigma) u)].\ee

Следовательно,
\be \Delta\varsigma '=div\Big[
\beta_1 (t,x,u,\varsigma) \nabla\varsigma + \fr{\beta_1}{x} (t,x,u,\varsigma) \varsigma +\fr{\beta_1}{u} (t,x,u,\varsigma) \varsigma\nabla u +\fr{\beta_1}{\varsigma} (t,x,u,\varsigma) \varsigma\nabla\varsigma $$ $$ + \gamma(t,x,u,\varsigma) \nabla u+\fr \gamma x (t,x,u,\varsigma) u +\fr \gamma u (t,x,u,\varsigma) u \nabla u +\fr \gamma \varsigma (t,x,u,\varsigma) u \nabla\varsigma\Big]. \ee

Обозначим $$v=u-\varphi, \tau=\varsigma - \psi,$$
$$\beta(t,x,v,\tau)=\fr{\beta_1}{u} (t,x,v+\varphi,\tau+\psi)
(\tau+\psi)+ \gamma(t,x,v+\varphi,\tau+\psi)$$ $$+\fr \gamma u
(t,x,v+\varphi,\tau+\psi) (v+\varphi), $$
$$\mu(t,x,v,\tau)=\beta_1 (t,x,v+\varphi,\tau+\psi)  +
\fr{\beta_1}{\varsigma} (t,x,v+\varphi,\tau+\psi) (\tau+\psi)
$$ $$+\fr \gamma \varsigma (t,x,v+\varphi,\tau+\psi) (v+\varphi),$$
$$g(t,x,v,\tau)=-\nabla \psi '+\fr{\beta_1}{u}
(t,x,v+\varphi,\tau+\psi) (\tau+\psi)\nabla \varphi+
\gamma(t,x,v+\varphi,\tau+\psi) \nabla \varphi$$ $$+\fr \gamma u
(t,x,v+\varphi,\tau+\psi) (v+\varphi) \nabla  \varphi
 +\beta_1 (t,x,v+\varphi,\tau+\psi) \nabla\psi $$ $$+ \fr{\beta_1}{\varsigma} (t,x,v+\varphi,\tau+\psi) (\tau+\psi)\nabla\psi +\fr \gamma \varsigma (t,x,v+\varphi,\tau+\psi) (v+\varphi) \nabla\psi $$ $$+ \fr{\beta_1}{x} (t,x,v+\varphi,\tau+\psi) (\tau+\psi)+\fr \gamma x (t,x,v+\varphi,\tau+\psi) (v+\varphi).$$

Тогда можно переписать (3.12) в следующем виде
:

\be \Delta\tau ' = div[\beta(t,x,v,\tau) \nabla v + \mu(t,x,v,\tau) \nabla \tau +g(t,x,v,\tau)].\ee

Аналогично, вводя обозначения $$D(t,x,v,\tau)=D_1(t,x,v+\varphi,\tau+\psi),$$ $$E(t,x,v,\tau)=E_1(t,x,v+\varphi,\tau+\psi),$$
$$f(t,x,v,\tau)=-\nabla \Delta^{-1}\varphi '+ D_1(t,x,v+\varphi,\tau+\psi) \nabla \varphi $$ $$+ E_1(t,x,v+\varphi,\tau+\psi) \nabla \psi+ (v+\varphi)M_1(t,x,v+\varphi,\tau+\psi),$$
перепишем (3.6) в форме \be v' = div[D(t,x,v,\tau) \nabla v + E(t,x,v,\tau) \nabla \tau +f(t,x,v,\tau)]. \ee

Заметим, что граничные и начальные условия на $v$ и $\tau$ получаются следующими:
\be v | _ {t=0} =v | _ {t=T},\ee \be \tau| _ {t=0} =
0, \ee \be v | _ {\partial \Omega} =0,\ \tau| _ {\partial \Omega} =
0.\ee

\textbf{Определение 3.1.} Пара функций $(v,\tau)$ из класса \be v \in W(\Omega,T), \tau\in
H^1(0,T;H^1_0(\Omega))\ee называется {\it слабым} решением задачи (3.13)-(3.17), если они удовлетворяют (3.15), (3.16), а равенства (3.13), (3.14) выполняются в пространстве $H^{-1}(\Omega)$ п.в. на
$(0,T)$.

Условия (3.15), (3.16) имеют смысл благодаря вложениям $$W\subset C([0,T];L_2(\Omega)),\ H^1(0,T;H^1_0(\Omega))\subset
C([0,T];H^1_0(\Omega)).$$ Условие (3.17) ''спрятано'' в пространстве $H^1_0(\Omega)$.

Как видим, всякое "хорошее"$,$ регулярное решение $(u,\varsigma)$
задачи (3.3)-(3.8) поро\-ждает пару $(v,\tau)$, которая, в
частности, является слабым решением задачи (3.13)-(3.17). Обратно,
пусть некая пара $(v,\tau)$ является слабым решением задачи
(3.13)-(3.17). Пусть функции $v,\tau,\psi,\varphi$ достаточно
регулярны, и имеют место соотношения (3.5), (3.9) и (3.10). Тогда
соответствующая пара $(u=v+\varphi,\varsigma=\tau+\psi)$
удовлетво\-ряет (3.6) и (3.11); последнее влечет (3.7) в силу
(3.9). Более того, имеют место (3.3), (3.4) и (3.8). А тогда пара
$(u,\sigma=\varsigma+\int\limits_0^u \nu_0(y)dy)$ удовлетворяет
(3.1) и (3.2).

При работе со слабыми решениями не обязательно считать, что все известные функции и область $\Omega$ регулярны. Достаточно условий, которые мы сейчас опишем.

Ради общности будем считать $\beta$ и $\mu$ матричными функциями.

Пусть $\Omega\subset\mathbb{R}^n$,
$n\in \mathbb{N}$, является произвольным ограниченным открытым множес\-твом, удовлетворяющим условию \be X(\Omega)\subset W^1_{p_0}(\Omega)\ee
с каким-нибудь $p_0>2$ (это условие выполнено при минимальных требованиях регулярности области).

Ниже мы предполагаем также, что
\\
i) $D$, $E$, $\mu$, $\beta: \R^{n+3}\to \R^{n\times n};$ $f,g: \R^{n+3}\to \R^n$.\\
ii) Каждая из этих шести функций (например, $D(t,x,v,\tau)$) измерима по
 $(t,x)$ при фиксированных $(v,\tau)$.\\
iii) Каждая из этих функций непрерывна по
 $(v,\tau)$ при фиксированных $(t,x)$.\\
iv) Имеют место оценки
\be |D(t,x,v,\tau)|\leq K_D,\ee
\be |E(t,x,v,\tau)|\leq K_E,\ee
\be |\beta(t,x,v,\tau)|\leq K_\beta,\ee
\be |\mu(t,x,v,\tau)|\leq K_\mu,\ee
\be |f(t,x,v,\tau)|\leq K_f|\tau|+ \widetilde f(t,x),\ee
\be |g(t,x,v,\tau)|\leq K_g(|v|+|\tau|)+ \widetilde g(t,x)\ee
с некоторыми константами $K_D,\dots,K_g$ и известными функциями $\widetilde f, \widetilde g \in L_2((0,T)\times\Omega)$.\\
v) \be (D(t,x,v,\tau)\xi, \xi)_{\R^n} \geq d(\xi, \xi)_{\R^n},\ee
где $d>0$ не зависит от $(t,x,v,\tau)\in \R^{n+3}$ и $\xi\in  \R^{n}$.

\begin{theorem} Найдется такое $T_0>0$, что при любом $T\leq T_0$ существует слабое решение задачи (3.13)-(3.17) в классе (3.18). \end{theorem}

Следующий пункт посвящен ее доказательству.

Так как краевая задача (3.13)-(3.14), (3.17) на отрезках любой (например, еди\-ничной) длины имеет решение, удовлетворяющее заданному условию $ v | _ {t=t_0} =v_{t_0}\in L_2,\ \tau| _ {t=t_0} =
\tau_{t_0}\in H^1_0$ в начале отрезка \cite{p2}, то полученное в теореме 3.1 решение можно, шаг за шагом, продлить на весь положительный луч, и мы получим

\begin{cor} Найдется такое $T_0>0$, что при любом $T\leq T_0$ существует пара функций \be v \in L_{2,loc}(0,\infty;H^1_0(\Omega))\bigcap H^1_{loc}(0,\infty;H^{-1}(\Omega)), \tau\in
H^1_{loc}(0,\infty;H^1_0(\Omega)),\ee которые удовлетворяют (3.15), (3.16), а также (3.13), (3.14) в $H^{-1}(\Omega)$ п.в. на
$(0,\infty)$. \end{cor}

Отметим также, что при моделировании реальных диффузионных процессов в полимерах всегда можно считать условия iv) выполненными\footnote{Естественность остальных наложенных нами условий ясна и не требует подробного обсуждения.}. В самом деле, по физическому смыслу \be 0\leq u(t,x)\leq 1\ee
(ведь $u$ это концентрация, так что она не менее $0\%$ и не более $100\%$). Кроме того, время релаксации в реальности положительно и ограничено, поэтому можно считать, что его обратное $\beta_0$ не менее некого $\beta_G>0$ (как во Введении), а потому $\beta_1 \leq -\beta_G$.
Равенства (3.7) и (3.8) дают представление
$$ \varsigma(t,x)=\varsigma_0(x)\exp\left(\int\limits_0^t \beta_1(\xi,x,u(\xi,x),\varsigma(\xi,x))\, d \xi\right)$$ $$+\int\limits_{0}^t
\exp\left(\int\limits_s^t \beta_1(\xi,x,u(\xi,x),\varsigma(\xi,x))\, d \xi\right) \gamma(s,x,u(s,x),\varsigma(s,x)) u(s,x)
\, ds.$$ Если $|\varsigma_0(x)|$ и $\gamma$ равномерно ограничены, $\varsigma$ также равномерно ограничено:
\be|\varsigma(t,x)|\leq e^{-t\beta_G}|\varsigma_0(x)|+C(\gamma)\int\limits_{0}^t
e^{(s-t)\beta_G}
\, ds\leq e^{-t\beta_G}|\varsigma_0(x)|+\frac {C(\gamma)}{\beta_G}.\ee
Поэтому коэффициенты системы (3.1)--(3.2) (и эквивалентной ей (3.6)--(3.7)) можно определить из эксперимента лишь на ограниченных $u$ и $\varsigma$, а "на бесконечности" мы их можем факти\-чески выбрать по своему усмотрению; например, ограниченными или даже имеющими доста\-точно быстро убыва\-ющие на бесконеч\-ности частные производные: так, чтобы выполнялись условия iv) на прео\-бразованные коэффи\-циенты.

Однако вопрос о том, при каких ограничениях на коэффициенты
условие (3.28), будучи заданным для граничной функции $\varphi$ и
для некого начального данного $u_0$, будет выполнено во все
последующие моменты времени для решения $u$ задачи (3.1)--(3.3) с
$ u | _ {t=0} =u_{0}$, с математической точки зрения изучен не до
конца (а ведь только такие решения, а потому и коэффициенты, имеют
физический смысл). Если $E_0$ не зависит\footnote{Как правило,
считают, что $E_0$ вообще зависит только от концентрации, см.
Введение.} от $\sigma$, а только от $t,x,u$, и $E_0(t,x,0)\equiv
0$, то \cite[теорема 7.4]{am2} решение сохраняет
неотрицательность. Сделав замену $\tilde u=1-u$,
$\tilde\sigma=-\sigma$ можно аналогичным образом доказать (при
предположениях $E_0(t,x,1)\equiv 0$ и $M_0(t,x,1,\sigma)=const$),
что условие $u\leq 1$ также сохраняется с течением времени. Вместе
с тем, нельзя наверняка утверждать, что другие $E_0$ не имеют
смысла. Поэтому мы рассматриваем модель в общем виде,
ограничиваясь лишь требо\-ваниями i)-v).

\section {Доказательство теоремы 3.1}

Нам понадобится один абстрактный результат. Пусть даны два
гильбертовых пространства $Z\subset Y$ с непрерывным вложением
$i:Z\to Y$, и $i(Z)$ плотно в $Y$. Стандартная схема
\cite{tem,gruy} отождествления $Y$ со своим сопряженным
пространством влечет
\begin{equation} Z\subset Y\equiv Y^*\subset Z^*,\end{equation}
где оба вложения плотны и непрерывны. При этом скалярное
произведение любых $f\in Y$ и $u\in Z$ в  $Y$ совпадает со
значением функционала $f$ из $Z^*$ на элементе $u\in Z$:
\begin{equation} (f,u)_Y=\langle f,u \rangle. \end{equation}

\begin{lemma} Пусть $Z$ сепарабельно, и дан непрерывный линейный оператор $A:Z\to
Z^*$, для которого найдется такое $\alpha>0$, что
\begin{equation} \langle Au,u\rangle\geq \alpha \|u\|^2_Z,\quad \forall u\in
Z.\end{equation} Тогда для любых $a\in Y$ и $f\in L_2(0,T;Z^*)$, $T>0$,
задача
\begin{gather}u'(t)+Au(t)=f(t),\, t\in(0,T),\\u(T)+a=u(0)\end{gather} имеет единственное решение в классе \begin{equation}u\in L_2(0,T;Z)\bigcap C([0,T];Y),\, u'\in
L_2(0,T;Z^*).\end{equation}\end{lemma}

\textbf{Доказательство.} Для любых $b\in Y$ и $f\in L_2(0,T;Z^*)$
задача Коши (4.4), $u(0)=b$ однозначно разрешима в классе (4.6) по теореме 1.1 из
\cite{ggz}, Глава VI, или по лемме 3.1.3 из \cite{gruy}. Обозначим
ее решение через $u_b$. Рассмотрим оператор $$\mathcal{G}: b
\mapsto u_b(T) + a.$$ Неподвижные точки этого оператора
соответствуют решениям задачи (4.4), (4.5).
Поэтому достаточно проверить, что $\mathcal{G}$ является сжатием на $Y$ (и применить принцип Банаха).

Пусть $b_1,b_2\in Y$. Обозначим $w=u_{b_1}-u_{b_2}$. Тогда $w'+Aw=0.$
Следовательно, для функции $\bar{w}(t)=e^{kt}w(t)$, $k>0$, выполнено равенство
\be \bar{w}'-k\bar{w}+A\bar{w}=0.\ee
Умножим (4.7) на $\bar{w}(t)$ при п.в. $t\in (0,T)$ (в смысле двойственности $Z^*$ и $Z$):
\be \langle\bar{w}', \bar{w}\rangle-k\langle\bar{w}, \bar{w}\rangle+\langle A\bar{w}, \bar{w}\rangle=0.\ee
При малом $k$, в силу (4.2), (4.3) и непрерывности вложения $i$, (4.8) влечет
\be \langle\bar{w}', \bar{w}\rangle\leq 0.\ee
Лемма III.1.2 из \cite{tem} дает, что функция $\|\bar{w}(t)\|^2_Y$ является невозрастающей.
Таким образом, $\|\bar{w}(T)\|_Y \leq \|\bar{w}(0)\|_Y$, т.е. $\|w(T)\|_Y\leq e^{-kT}\|w(0)\|_Y$ и $\mathcal{G}$ является сжатием.
$ \Box $

Рассмотрим теперь следующую вспомогательную задачу:

\be \fr v t + \varepsilon \Delta^2 v= \lambda div[D(t,x,v,\tau) \nabla v + E(t,x,v,\tau) \nabla \tau +f(t,x,v,\tau)], \ee \be \fr
{\tau} t + \varepsilon \Delta^2 \tau= \lambda \Delta^{-1}div[\beta(t,x,v,\tau) \nabla v + \mu(t,x,v,\tau) \nabla \tau +g(t,x,v,\tau)], \ee \be v | _ {t=0} =v | _ {t=T},\ee \be \tau| _ {t=0} =0.\ee Тут $\varepsilon>0$, $\lambda\in[0,1]$ --- параметры.

\textbf{Определение 4.1.} Пара функций $(v,\tau)$ из класса \be v \in W_2(\Omega,T),\tau\in
W_1(\Omega,T)\ee называется {\it слабым} решением задачи (4.10)-(4.13), если равенство
(4.10) имеет место в пространстве $H^{-2}(\Omega)$ п.в. на
$(0,T)$, (4.11) выполнено в $H^{-1}(\Omega)$ п.в. на
$(0,T)$, (4.12) верно в $L_2(\Omega)$, а (4.13) --- в $H^1_0(\Omega)$.

Последние два условия имеют смысл, т.к. $$W_1\subset C([0,T];H^1_0(\Omega)),\ W_2\subset
C([0,T];L_2(\Omega)).$$

\begin{lemma} Найдется такое $T_0>0$, что при любом $T\leq T_0$ для всякого слабого решения $(v,\tau)$ задачи
(4.10)-(4.13) выполнена априорная оценка: \be
\varepsilon\|v\|^2_{L_2(0,T;H^2_0(\Omega))}+ \varepsilon\|\tau\|^2_{L_2(0,T;X)}+$$ $$+ \lambda \|v\|^2_{L_2(0,T;H^1_0(\Omega))} +
\|\tau\|^2_{L_2(0,T;H^1_0(\Omega))}\leq C, \ee где $C$ не
зависит от $\lambda$ и $\varepsilon$. \end{lemma}
\textbf{Доказательство.} Рассмотрим произведения слагаемых (4.11) с $-\Delta\tau(t)\in H_0^1(\Omega)$ при п.в.
$t\in[0,T]$  (в смысле двойственности
 $H^{-1}(\Omega)$ и $H_0^1(\Omega)$): \be -\left\langle  \tau', \Delta\tau\right\rangle   - \left\langle \varepsilon \Delta^2 \tau, \Delta\tau\right\rangle $$ $$ =- \lambda \left( \Delta^{-1}div[\beta(t,x,v,\tau) \nabla v + \mu(t,x,v,\tau) \nabla \tau +g(t,x,v,\tau)], \Delta\tau\right) . \ee
Но \be-\left\langle  \tau', \Delta\tau\right\rangle=\left\langle  \tau',\tau\right\rangle_1=\frac 1 2\frac {d}{dt}\|\tau\|^2_1\ee (см. \cite[Лемма  III.1.2]{tem}).
Таким образом,
\be \frac 1 2\frac {d}{dt}\|\tau\|^2_1   + \varepsilon(\nabla\Delta \tau, \nabla\Delta\tau) $$ $$ =- \lambda \left\langle div[\beta(t,x,v,\tau) \nabla v + \mu(t,x,v,\tau) \nabla \tau +g(t,x,v,\tau)], \tau\right\rangle .\ee
Обозначим $\bar{v}(t)=e^{-kt}v(t)$, $\bar{\tau}(t)=e^{-kt}\tau(t)$, где $k>0$ будет определено позднее. Тогда
\be \frac 1 2\frac {d}{dt}\|e^{kt}\bar{\tau}\|^2_1   + e^{2kt}\varepsilon(\nabla\Delta \bar\tau, \nabla\Delta\bar\tau) $$ $$ =\lambda \Big( \beta(t,x,e^{kt}\bar v(t),e^{kt}\bar\tau(t)) \nabla \bar ve^{kt} + \mu(t,x,e^{kt}\bar v(t),e^{kt}\bar \tau(t)) \nabla \bar\tau e^{kt} $$ $$+g(t,x,e^{kt}\bar v(t),e^{kt}\bar\tau(t)), \nabla\bar \tau(t)e^{kt}\Big) .\ee
Теперь обозначим $$\beta_k(t,x,\bar v(t),\bar\tau(t))=\beta(t,x,e^{kt}\bar v(t),e^{kt}\bar\tau(t)),$$
$$\mu_k(t,x,\bar v(t),\bar\tau(t))=\mu(t,x,e^{kt}\bar v(t),e^{kt}\bar\tau(t)),$$ $$g_k(t,x,\bar v(t),\bar\tau(t))= e^{-kt}g(t,x,e^{kt}\bar v(t),e^{kt}\bar\tau(t)).$$
Эти функции удовлетворяют оценкам типа (3.22), (3.23), (3.25) с
теми же константами.

Имеем:
\be \frac 1 2\frac {d}{dt}\|\bar{\tau}\|^2_1 + k \|\bar{\tau}\|^2_1    + \varepsilon(\bar\tau, \bar\tau)_X $$ $$ = \lambda\Big(\beta_k(t,x,\bar v(t),\bar\tau(t)) \nabla \bar v + \mu_k(t,x,\bar v(t),\bar \tau(t)) \nabla \bar\tau $$ $$+g_k(t,x,\bar v(t),\bar\tau(t)), \nabla\bar \tau(t)\Big) .\ee
Интегрируя по интервалу $(0,T)$, получим
\be   \frac 1 2\|\bar\tau(T)\|^2_1  +k\int\limits_0^T \|\bar{\tau}(s)\|^2_1\,ds+ \varepsilon\int\limits_0^T\|\bar\tau(s)\|^2_X\,ds  $$ $$ = \lambda\int\limits_0^T\Big( \beta_k(s,x,\bar v(s),\bar \tau(s)) \nabla \bar v(s) + \mu_k(s,x,\bar v(s),\bar\tau(s)) \nabla \bar\tau(s) $$ $$+g_k(s,x,\bar v(s),\bar \tau(s)), \nabla\bar\tau(s)\Big)\,ds.  \ee
Применяя неравенство Коши-Буняковского и неравенство Коши $ab\leq c a^2+ \frac 1 {4c} b^2$ и отбрасывая первое слагаемое, придем к неравенству
\be   k\int\limits_0^T \|\bar{\tau}(s)\|^2_1\,ds+ \varepsilon\int\limits_0^T\|\bar\tau(s)\|^2_X\,ds  $$ $$ \leq \frac{ \lambda K^2_\beta}4\int\limits_0^T \|\bar{v}(s)\|^2_1\,ds + \lambda \int\limits_0^T \|\bar{\tau}(s)\|^2_1\,ds+\lambda K_\mu\int\limits_0^T \|\bar{\tau}(s)\|^2_1\,ds$$ $$+\frac {\lambda} 4 \int\limits_0^T \|g_k(s,\cdot,\bar v(s,\cdot),\bar \tau(s,\cdot))\|^2\,ds+ \lambda\int\limits_0^T \|\bar{\tau}(s)\|^2_1\,ds.  \ee
Заметим, что $$\int\limits_0^T \|g_k(s,\cdot,\bar v(s,\cdot),\bar \tau(s,\cdot))\|^2\,ds\leq  \int\limits_0^T \|K_g[|\bar v(s,\cdot)|+|\bar \tau(s,\cdot)|]+\widetilde g(s,\cdot)\|^2\,ds$$ $$\leq 3 K_g^2\int\limits_0^T \|\bar v(s,\cdot)\|^2\,ds + 3 K_g^2\int\limits_0^T \|\bar \tau(s,\cdot)\|^2\,ds + 3\int\limits_0^T \|\widetilde g(s,\cdot)\|^2\,ds$$ $$\leq 3 K_g^2 K_\Omega^2\int\limits_0^T \|\bar v(s)\|^2_1\,ds + 3 K_g^2 K_\Omega^2 \int\limits_0^T \|\bar \tau(s)\|_1^2\,ds + 3\|\widetilde g\|^2_{L_2((0,T)\times \Omega)}.$$
Поэтому\be   (k-2-K_\mu-\frac 3 4 K_g^2 K_\Omega^2)\int\limits_0^T \|\bar{\tau}(s)\|^2_1\,ds+ \varepsilon\int\limits_0^T\|\bar\tau(s)\|^2_X\,ds  $$ $$ \leq\lambda (\frac {K^2_\beta} 4+\frac 3 4 K_g^2 K_\Omega^2)\int\limits_0^T \|\bar{v}(s)\|^2_1\,ds +\frac 34\|\widetilde g\|^2_{L_2((0,T)\times \Omega)}.  \ee

Возьмем $k\geq 4+2K_\mu+\frac 3 2 K_g^2 K_\Omega^2$.

Тогда (4.23) влечет неравенство
\be   \int\limits_0^T \|\bar{\tau}(s)\|^2_1\,ds\leq \frac {C_1} k(1+\lambda\int\limits_0^T \|\bar{v}(s)\|^2_1\,ds), \ C_1=\frac {K^2_\beta} 2+\frac 3 2 K_g^2 K_\Omega^2 + \frac 32\|\widetilde g\|^2_{L_2((0,T)\times \Omega)}.  \ee

Положим $T_0=\frac 1 k$. Тогда $\|\tau(t)\|_1\leq e \|\bar{\tau} (t)\|_1$ при $0\leq t \leq T\leq T_0$, и из (4.24) следует, что

\be   \int\limits_0^T \|\tau(s)\|^2_1\,ds\leq \frac {C_1 e} k(1+\lambda\int\limits_0^T \|v(s)\|^2_1\,ds).  \ee

Рассмотрим теперь произведения слагаемых (4.10) с $v(t)\in H_0^2(\Omega)$ при п.в.
$t\in[0,T]$  (в смысле двойственности
 $H^{-2}(\Omega)$ и $H_0^2(\Omega)$): \be \left\langle  v', v\right\rangle   + \left\langle \varepsilon \Delta^2 v, v\right\rangle $$ $$ = \lambda \left\langle div[D(t,x,v,\tau) \nabla v + E(t,x,v,\tau) \nabla \tau +f(t,x,v,\tau)], v\right\rangle . \ee
Опять, в силу \cite[Лемма  III.1.2]{tem}, мы имеем \be \left\langle  v',v\right\rangle=\frac 1 2\frac {d}{dt}\|v\|^2,\ee
и поэтому
\be \frac 1 2\frac {d}{dt}\|v\|^2   + \varepsilon(\Delta v, \Delta v) $$ $$ =-\lambda (D(t,x,v,\tau) \nabla v + E(t,x,v,\tau) \nabla \tau +f(t,x,v,\tau), \nabla v) .\ee

Следовательно, \be    \varepsilon\int\limits_0^T\|v(s)\|^2_2\,ds
$$ $$ = - \lambda \int\limits_0^T\Big( D(s,x,v(s),\tau(s)) \nabla
v(s) + E(s,x,v(s),\tau(s)) \nabla \tau(s) $$ $$+f(s,x, v(s),
\tau(s)), \nabla v(s)\Big)\,ds.  \ee С помощью неравенства Коши,
(3.21) и (3.26), получаем неравенство \be
\varepsilon\int\limits_0^T\|v(s)\|^2_2\,ds  + \lambda d
\int\limits_0^T (\nabla  v(s), \nabla v(s))\,ds$$ $$\leq \frac
{3\lambda K^2_E}{4d}\int\limits_0^T \|{\tau}(s)\|^2_1\,ds +
\frac{\lambda d}{3}\int\limits_0^T \|{v}(s)\|^2_1\,ds$$ $$+\frac
{3\lambda} {4d} \int\limits_0^T
\|f(s,\cdot,v(s,\cdot),\tau(s,\cdot))\|^2\,ds+  \frac{\lambda
d}{3}\int\limits_0^T \|{v}(s)\|^2_1\,ds.  \ee По тому же принципу,
как для $g_k$ выше, мы имеем $$\int\limits_0^T \|f(s,\cdot,
v(s,\cdot), \tau(s,\cdot))\|^2\,ds$$ $$  \leq  3 K_f^2 K_\Omega^2
\int\limits_0^T \|\tau(s)\|_1^2\,ds + 3\|\widetilde
f\|^2_{L_2((0,T)\times \Omega)}.$$ Итак, из (4.30) и (4.25)
заключаем \be    \varepsilon\int\limits_0^T\| v(s)\|^2_2\,ds  +
\frac{\lambda d}{3} \int\limits_0^T \| v(s)\|_1^2\,ds$$ $$\leq
(\frac {3K^2_E}{4d}+ \frac 9{4d}K_f^2 K_\Omega^2)\int\limits_0^T
\|{\tau}(s)\|^2_1\,ds +\frac 9 {4d} \|\widetilde
f\|^2_{L_2((0,T)\times \Omega)}.  $$ $$ \leq \frac {C_2}
k(1+\lambda\int\limits_0^T \|{v}(s)\|^2_1\,ds)+C, \ee $$ C_2=C_1
e(\frac {3K^2_E}{4d}+ \frac 9{4d}K_f^2 K_\Omega^2).$$

Без ограничения общности можно считать, что $k\geq \frac {6C_2}d$.  Тогда (4.31) влечет $$\frac{\lambda d} 6\int\limits_0^T \|\bar{v}(s)\|^2_1\,ds \leq C.$$

Таким образом, правые части (4.23) и (4.31) ограничены, и мы приходим к (4.15).
$ \Box $

\begin{lemma} Пусть слабое решение $(v,\tau)$ задачи
(4.10)-(4.13) удовлетворяет оценке (4.15). Тогда имеет место следующая оценка производных по времени: \be
 \|v'\|_{L_2(0,T;H^{-2}(\Omega))} +
\|\tau'\|_{L_2(0,T;H^{-1}(\Omega))}\leq C(1+\sqrt{\varepsilon}) \ee где $C$ не
зависит от $\lambda$ и $\varepsilon$. \end{lemma}
\textbf{Доказательство.} Действительно, так как $H^{-1}(\Omega)\subset H^{-2}(\Omega)$ непрерывным образом, (4.10) и (4.15) влекут $$\|v'\|_{L_2(0,T;H^{-2}(\Omega))}\leq \varepsilon \|\Delta^2 v\|_{L_2(0,T;H^{-2}(\Omega))}+$$ $$ \lambda C \|div[D(t,x,v,\tau) \nabla v + E(t,x,v,\tau) \nabla \tau +f(t,x,v,\tau)]\|_{L_2(0,T;H^{-1}(\Omega))}$$
$$\leq \sqrt\varepsilon \sqrt\varepsilon\| v\|_{L_2(0,T;H^{2}_0(\Omega))}+$$ $$ \sqrt\lambda C \|D(t,x,v,\tau) \nabla v + E(t,x,v,\tau) \nabla \tau +f(t,x,v,\tau)\|_{L_2(0,T;L_2(\Omega))}$$ $$\leq C\sqrt\varepsilon + C \sqrt\lambda[K_D \| v\|_{L_2(0,T;H^1_0(\Omega))} + K_E\|\tau\|_{L_2(0,T;H^1_0(\Omega))} +\|f(t,x,v,\tau)\|_{L_2(0,T;L_2(\Omega))}]$$ $$\leq C\sqrt\varepsilon + C\sqrt\lambda [K_D \| v\|_{L_2(0,T;H^1_0(\Omega))} + K_E\|\tau\|_{L_2(0,T;H^1_0(\Omega))} $$ $$+K_f\|\tau\|_{L_2(0,T;L_2(\Omega))}+\|\widetilde f\|_{L_2((0,T)\times\Omega))}]$$ $$\leq C\sqrt\varepsilon + C\sqrt\lambda [ \| v\|_{L_2(0,T;H^1_0(\Omega))} +\|\tau\|_{L_2(0,T;H^1_0(\Omega))} +1]$$ $$\leq C(1+\sqrt\varepsilon).$$

Аналогично, раз $H^{1}_0(\Omega)\subset H^{-1}(\Omega)$, (4.11) и (4.15) дают
$$\|\tau'\|_{L_2(0,T;H^{-1}(\Omega))}\leq \varepsilon \|\Delta^2
\tau\|_{L_2(0,T;H^{-1}(\Omega))}+$$ $$ \lambda C
\|\Delta^{-1}div[\beta(t,x,v,\tau) \nabla v + \mu(t,x,v,\tau)
\nabla \tau +g(t,x,v,\tau)]\|_{L_2(0,T;H^{1}_0(\Omega))}$$
$$\leq \sqrt\varepsilon \sqrt\varepsilon\| \tau\|_{L_2(0,T;X)}+$$ $$C\sqrt\lambda \|\beta(t,x,v,\tau) \nabla v + \mu(t,x,v,\tau) \nabla \tau +g(t,x,v,\tau)\|_{L_2(0,T;L_2(\Omega))}$$ $$\leq C\sqrt\varepsilon + C\sqrt\lambda [K_\beta \| v\|_{L_2(0,T;H^1_0(\Omega))} + K_\mu\|\tau\|_{L_2(0,T;H^1_0(\Omega))} +\|g(t,x,v,\tau)\|_{L_2(0,T;L_2(\Omega))}]$$ $$\leq C\sqrt\varepsilon + C \sqrt\lambda[K_\beta \| v\|_{L_2(0,T;H^1_0(\Omega))} + K_\mu\|\tau\|_{L_2(0,T;H^1_0(\Omega))} $$ $$+K_g\|v\|_{L_2(0,T;L_2(\Omega))}+K_g\|\tau\|_{L_2(0,T;L_2(\Omega))}+\|\widetilde g\|_{L_2((0,T)\times\Omega))}]$$ $$\leq C\sqrt\varepsilon + C\sqrt\lambda [ \| v\|_{L_2(0,T;H^1_0(\Omega))} +\|\tau\|_{L_2(0,T;H^1_0(\Omega))} +1]$$ $$\leq C(1+\sqrt\varepsilon).$$$ \Box $

\begin{lemma} При любом $T\leq T_0$ имеется слабое решение $(v,\tau)$ задачи
(4.10)-(4.13) в классе (4.14). \end{lemma}

 \textbf{Доказательство.} Введем вспомогательные операторы по формулам:
$$Q_1:W_2\times W_1\to L_2 (0, T; H^{-2}(\Omega)), $$
$$Q_1 (v, \tau) = div[D(\cdot,\cdot,v,\tau) \nabla v],$$
$$Q_2:W_2\times W_1\to L_2 (0, T; H^{-2}(\Omega)), $$
$$Q_2 (v, \tau) = div[E(\cdot,\cdot,v,\tau) \nabla \tau],$$
$$Q_3:W_2\times W_1\to L_2 (0, T; H^{-2}(\Omega)), $$
$$Q_3 (v, \tau) = div[f(\cdot,\cdot,v,\tau)],$$
$$Q_4:W_2\times W_1\to L_2 (0, T; H^{-1}(\Omega)), $$
$$Q_4 (v, \tau) = \Delta^{-1}div[\beta(\cdot,\cdot,v,\tau) \nabla v],$$
$$Q_5:W_2\times W_1\to L_2 (0, T; H^{-1}(\Omega)), $$
$$Q_5 (v, \tau) = \Delta^{-1}div[\mu(\cdot,\cdot,v,\tau) \nabla \tau],$$
$$Q_6:W_2\times W_1\to L_2 (0, T; H^{-1}(\Omega)), $$
$$Q_6 (v, \tau) = \Delta^{-1}div[g(\cdot,\cdot,v,\tau)],$$
$$Q:W_2\times W_1\to L_2 (0, T; H^{-2}(\Omega))\times L_2 (0, T; H ^ {-1}(\Omega))\times L_2(\Omega)\times H^1_0(\Omega), $$
$$Q (v, \tau) = (-Q_1 (v,\tau) -Q_2(v, \tau) -Q_3(v,\tau), -Q_4 (v,\tau) -Q_5(v, \tau) -Q_6(v,\tau), 0,0), $$

$$\tilde{A}_1:W_1\to L_2 (0, T; H^{-1}(\Omega))\times  H^1_0(\Omega),$$ $$
\tilde{A}_1 (u)= (u'+\varepsilon \Delta^2 u, u | _ {t=0})
,
$$
$$\tilde{A}_2:W_2\to L_2 (0, T; H^{-2}(\Omega))\times  L_2(\Omega),$$ $$
\tilde{A}_2 (u)= (u'+\varepsilon \Delta^2
u, u | _ {t=0} - u| _{t=T}),
$$
$$\tilde{A}:W_2\times W_1\to L_2 (0, T; H^{-2}(\Omega))\times L_2 (0, T; H ^ {-1}(\Omega))\times L_2(\Omega)\times H^1_0(\Omega),$$ $$
\tilde{A} (v,\tau)= (v'+\varepsilon \Delta^2
v, \tau'+\varepsilon \Delta^2
\tau , v | _ {t=0} - v | _{t=T}, \tau | _ {t=0}).$$

Слабая постановка задачи (4.10) - (4.13) эквивалентна операторному уравнению
\begin{equation}\tilde {A} (v, \tau) + \lambda Q (v, \tau) = (0, 0, 0, 0) \end{equation}

При условии (3.19) оператор $Q $ вполне непрерывен (см.
\cite{p2}).

Заметим, что $$\left\langle \Delta^2 u, u\right\rangle =(\Delta u, \Delta u)=\|u\|_2^2$$ для $u\in H^2_0(\Omega)$, и
$$\left\langle \Delta^2 u, u\right\rangle_1=-\left\langle \Delta^2 u, \Delta u\right\rangle=(\nabla\Delta u, \nabla\Delta u)=\|u\|_X^2$$ для $u\in X$.
Поэтому оператор $
\tilde {A}_1 $ непрерывно обратим (по теореме 1.1 из
\cite{ggz}, Глава VI, или по лемме 3.1.3 из \cite{gruy}). А оператор $
\tilde {A}_2 $ непрерывно обратим по доказанной выше лемме 4.1. Поэтому $
\tilde {A}$ также обратим.

Перепишем (4.33) в виде

\begin{equation} (u, \tau) + \lambda\tilde {A} ^ {-1} Q (u, \tau) = (0, 0) \end{equation} Априорные оценки лемм 4.2 и 4.3 не позволяют уравнению (4.34) иметь
решения на границе достаточно большого шара $B $ в
$W_2\times W_1, $ не зависящего от $ \lambda. $  Рассмотрим степень Лере-Шаудера (см.
напр. \cite{kras,kz,ll}) вполне непрерывного векторного поля $I + \lambda \tilde {A} ^ {-1} Q $ ($I$ -- тождественное отображение) на шаре $B $ ---
$$deg _ {LS} (I + \lambda \tilde {A} ^ {-1} Q, B, 0). $$  В силу гомотопической инвариантности степени
$$ deg _ {LS} (I + \lambda \tilde {A} ^ {-1} Q, B, 0) = deg _ {LS} (I, B,
0) =1\neq 0. $$ Итак, уравнение (4.34) (а потому и задача
(4.10) - (4.13)) имеет решение в шаре $B $ при всяком $
\lambda $ (ср. с подобными рассуждениями в \cite{gruy,zvdkn}). $ \Box $

\textbf{Доказательство теоремы 3.1.}

Возьмем какую-нибудь убывающую последовательность
положительных чисел $\varepsilon_m\to 0 $. По лемме 4.4 найдутся  пары $(v_m, \tau_m) $, которые являются слабым решением задачи
(4.10) - (4.13) с $\lambda=1$,
$\varepsilon =\varepsilon_m$.
Заметим, что из (4.11) следует, что  \be \Delta\tau_m ' + \varepsilon_m \Delta^3 \tau_m = div [\beta (t, x, v_m, \tau_m) \nabla v_m + \mu (t, x, v_m, \tau_m) \nabla \tau_m +g (t, x, v_m, \tau_m)] \ee в $H ^ {-3} (\Omega) $ п.в. на
$ (0, T) $.
Благодаря априорной оценке
(4.15), без потери общности (переходя к подпоследовательности, если
необходимо) можно предположить, что существуют
пределы \\
$v_* = \lim\limits _ {m\to
\infty} ^ {} {v} _ {m} $,   \\
$\tau_* = \lim\limits _ {m\to \infty} ^ {} {\tau} _ {m} $, которые
являются  слабыми в $L _ {2} (0, T; H_0^1 (\Omega)) $. Кроме того,
в силу леммы 4.3, без ограничения общности можно считать, что $
{v} _m '\to v_* '$ слабо в $L_2 (0, T; H ^ {-2}) $, $ {\tau} _m '
\to \tau_* ' $ слабо в $L_2 (0, T; H ^ {-1}) $. Следовательно, $
{v} _m \to v_* $ слабо в $C ([0, T]; H ^ {-2}) $, $ {\tau} _m  \to
\tau_* $ слабо в $C ([0, T]; H ^ {-1}) $. Тогда, \cite [Следствие
4] {sim1}, $v_m\to v_ * $ сильно в $C ([0, T]; H ^ {-3}) $,
$\tau_m \to \tau_*$ сильно в $C ([0, T]; H ^ {-2}) $. Поэтому
$v_*$ и $\tau _*$ удовлетворяют (3.15), (3.16). Кроме того, \cite
[Следствие 4] {sim1}, $v_m\to v_ *, $ $\tau_m \to \tau_*$ сильно в
$L_2 (0, T; L_2) $. По теореме Красносельского \cite {kras, skry}
о непрерывности оператора суперпозиции
$$D (\cdot, \cdot, v_m, \tau_m) \to D (\cdot, \cdot, v_ *,\tau_ *), $$
$$E (\cdot, \cdot, v_m, \tau_m) \to E (\cdot, \cdot, v_ *,\tau_ *), $$
$$\beta (\cdot, \cdot, v_m, \tau_m) \to \beta (\cdot, \cdot, v_ *,\tau_ *), $$
$$\mu (\cdot, \cdot, v_m, \tau_m) \to \mu (\cdot, \cdot, v_ *,\tau_ *), $$
сильно в $L_p ((0, T) \times \Omega) ^ {n\times n} $ для любого $p <\infty$;
$$f (\cdot, \cdot, v_m, \tau_m) \to f (\cdot, \cdot, v_ *,\tau_ *), $$
$$g (\cdot, \cdot, v_m, \tau_m) \to г (\cdot, \cdot, v_ *,\tau_ *) $$
сильно в $L_2 ((0, T) \times \Omega) ^n$.

Если последовательность функций $y_m$ сходится слабо в $L_2 ((0, T) \times \Omega)$, и другая последовательность $z_m$ сходится сильно в $L_p ((0, T) \times \Omega) $, то их поточечные произ\-ведения $y_m z_m$ сходятся слабо к произведению их пределов в $L_q ((0, T) \times \Omega) $, $\frac 1 p + \frac 1 2 =\frac 1 q $.
Поэтому $$D (\cdot, \cdot, v_m, \tau_m) \nabla v_m\to D (\cdot, \cdot, v_ *,\tau_ *)\nabla v_ *, $$
$$E (\cdot, \cdot, v_m, \tau_m) \nabla \tau_m\to E (\cdot, \cdot, v_ *,\tau_ *)\nabla \tau_ *, $$
$$\beta (\cdot, \cdot, v_m, \tau_m) \nabla v_m\to \beta (\cdot, \cdot, v_ *,\tau_ *)\nabla v_ *, $$
$$\mu (\cdot, \cdot, v_m, \tau_m) \nabla \tau_m\to \mu (\cdot, \cdot, v_ *,\tau_ *)\nabla \tau_*$$ слабо в $L_q (0, T; L_q (\Omega) ^n) $ для $1\leq q <2 $.

Переходя к пределу при $m \to \infty$ в (4.10) с $\lambda=1$,
$\varepsilon =\varepsilon_m$, $v=v_m$, $\tau =\tau_m$ и в (4.35)
(например, в смысле распределений на $ (0, T) $ со значениями в $W
^ {-6} _q$), мы заключаем, что $(v_ *,\tau_ *) $ есть слабое
решение задачи (3.13)-(3.17). Остается заметить, что правые части
(и, следовательно, левые части) (3.13), (3.14) принадле\-жат $L_2
(0, T; H ^ {-1}) $. Тогда $v_* \in W (\Omega, T) $. Так как
$\tau_* ' \in H ^ {-1} (0, T; H ^ {1} _0), $ и $ \Delta\tau_* '
\in L_2 (0, T; H ^ {-1}) $, мы имеем $\tau_* ' \in L_2 (0, T; H ^
{1} _0) $. $\Box$

\section{Задача о периодических решениях}

В этом пункте исследуется периодическая задача для системы (3.1)-(3.3), то есть "временные" условия берутся в виде \be u(0,x) =u (T,x), \sigma(0,x)= \sigma(T,x), \ x\in\Omega\ee
с некоторым $T>0$. Предполагается, что функции $\beta_0, M_0, D_0, E_0$ периодичны с периодом $T$ по переменной $t$. Как и в пункте 3, после замены переменной $\sigma$ на $\varsigma$ мы приходим к системе (3.6), (3.7). Будем предполагать, что граничное условие на $\varsigma$ известно: \be \varsigma(t,x)=\psi(t,x),\ (t,x)\in [0,T]\times \partial\Omega, \ee и обе функции $\varphi$ и $\psi$ определены и периодичны с периодом $T$ по переменной $t$ на $[0,T]\times\overline{\Omega}$; при этом предполагается условие согласования (3.9) (на самом деле $\psi$ можно определить, зная $\varphi$, как периодическое решение уравнения (3.9) при фиксированном $x$, вопросы периодической разрешимости таких уравнений см. напр. в \cite{dem,kz}, но может быть проблема с единственностью такого решения и зависимостью от $x$; чтобы не заостряться на этом, мы считаем $\psi$ известным). После второй замены переменных и тех же преобразований, что и в пункте 3, мы приходим к системе (3.13), (3.14). При этом имеем граничное условие (3.17) и условие периодичности \be v | _ {t=0} =v | _ {t=T},\ \tau| _ {t=0} =
\tau | _ {t=T}. \ee

Определение слабого решения задачи  (3.13), (3.14), (3.17), (5.3) в классе (3.18) аналогично определению 3.1. При работе со слабыми решениями будем предпо\-лагать выполненными условия
пункта 3, за исключением (3.24) -- (3.26).
Вместо них ниже мы предполагаем более жесткие условия, а именно

v)
\be |f(t,x,v,\tau)|\leq  \widetilde f(t,x),\ |g(t,x,v,\tau)|\leq \widetilde g(t,x)\ee
с некоторыми известными функциями $\widetilde f, \widetilde g \in L_2((0,T)\times\Omega)$, а также

vi) существуют положительные числа $\Gamma$ и $\Gamma_0$ такие, что \be (D(\cdot)\xi,\xi)_{\R^n}- (\mu(\cdot)\eta,\eta)_{\R^n}+ \left(\left[E(\cdot)\Gamma-\frac {\beta(\cdot)}{\Gamma }\right]\xi,\eta\right)_{\R^n}\geq \Gamma_0(|\xi|^2+|\eta|^2)\ee
для любых $\xi,\eta\in \R^n$.

Добиться выполнения условия v) можно, применив рассуждения из конца пункта 3. Условие vi) уже
существенно ограничивает общность модели. Тем не менее, следующие неформальные аргументы показывают, что оно может быть выполнено для реаль\-ных ситуаций, причем надо брать достаточно большое $\Gamma$. Во первых (см. п. 3), $\mu$ это по сути скаляр, причем, считая $\gamma$ малым, можно считать его близким к $-\beta_0$, то есть $-\mu\approx \beta_0$. При больших концентрациях $\beta_0$ и $D$ велики (см. Введение), и можно считать (5.5) выполненным. А при малых концентрациях уже мало $E$, что при большом $\Gamma$ опять играет в пользу (5.5).

\begin{theorem} В описанных условиях при любом $T>0$ существует слабое решение задачи (3.13), (3.14), (3.17), (5.3) в классе (3.18). \end{theorem}

При доказательстве используются вспомогательная система (4.10)-(4.11) и следу\-ющая априорная оценка:

\begin{lemma} Для всякого решения $(v,\tau)$ задачи
(4.10), (4.11), (5.3) в классе (4.14) выполнена априорная оценка: \be
\varepsilon\|v\|^2_{L_2(0,T;H^2_0(\Omega))}+ \varepsilon\|\tau\|^2_{L_2(0,T;X)}+$$ $$+ \lambda \|v\|^2_{L_2(0,T;H^1_0(\Omega))} + \lambda
\|\tau\|^2_{L_2(0,T;H^1_0(\Omega))}\leq C, \ee где $C$ не
зависит от $\lambda$ и $\varepsilon$. \end{lemma}

\textbf{Доказательство.} Условие (5.5) можно переписать в виде \be (D(\cdot)\Gamma^2\xi,\xi)_{\R^n}- (\mu(\cdot)\eta,\eta)_{\R^n}+ \left(E(\cdot)\Gamma^2\xi, \eta\right)_{\R^n}-({\beta(\cdot)}\xi,\eta)_{\R^n}$$ $$\geq \Gamma_0(\Gamma^2|\xi|^2+|\eta|^2)\ee
для любых $\xi,\eta\in \R^n$ (для этого достаточно подставить в (5.5) $\Gamma\xi$ вместо $\xi$).

Как в доказательстве леммы 4.2, мы имеем (4.28). Кроме того, мы имеем (4.18), и, следовательно, \be \frac 1 2\frac {d}{dt}\|\tau\|^2_1    + \varepsilon(\tau, \tau)_X $$ $$ =\lambda\Big(\beta(t,x,v,\tau) \nabla v + \mu(t,x, v, \tau) \nabla \tau, \nabla \tau\Big) $$ $$+ \lambda\Big(g(t,x, v,\tau), \nabla \tau\Big).\ee

Сложим это с (4.28), умноженным на $\Gamma^2$:
 \be \frac {\Gamma^2} 2\frac {d}{dt}\|v\|^2   + \Gamma^2\varepsilon(\Delta v, \Delta v)+ \frac 1 2\frac {d}{dt}\|\tau\|^2_1    + \varepsilon(\tau, \tau)_X $$ $$ =-\lambda (D(t,x,v,\tau)\Gamma^2 \nabla v + E(t,x,v,\tau) \Gamma^2\nabla \tau, \nabla v) $$ $$+ \lambda\Big(\beta(t,x,v,\tau) \nabla v + \mu(t,x, v, \tau) \nabla \tau, \nabla \tau\Big) $$ $$+ \lambda\Big(g(t,x, v,\tau), \nabla \tau\Big)-\lambda\Gamma^2(f(t,x,v,\tau), \nabla v)  .\ee

Используя (5.7), заключаем, что

\be \frac {\Gamma^2} 2\frac {d}{dt}\|v\|^2   + \Gamma^2\varepsilon(\Delta v, \Delta v)+ \frac 1 2\frac {d}{dt}\|\tau\|^2_1    + \varepsilon(\tau, \tau)_X $$ $$+\lambda\Gamma_0(\Gamma^2\|\nabla v\|^2+\|\nabla \tau\|^2) $$ $$ \leq \lambda\Big(g(t,x, v,\tau), \nabla \tau\Big)-\lambda\Gamma^2(f(t,x,v,\tau), \nabla v)  .\ee

Интегрируя по интервалу $(0,T)$, получим \be
 \Gamma^2\varepsilon\int\limits_0^T \|v(s)\|^2_2\,ds+ \varepsilon\int\limits_0^T\|\tau(s)\|^2_X\,ds + \lambda\Gamma_0\Gamma^2\int\limits_0^T \|\nabla v(s)\|^2\,ds+ \lambda\Gamma_0\int\limits_0^T\|\nabla\tau(s)\|^2\,ds $$ $$ \leq \lambda\int\limits_0^T\Big(g(s,x,v(s),\tau(s)), \nabla\tau(s)\Big)\,ds - \lambda\Gamma^2\int\limits_0^T\Big(f(s,x,v(s),\tau(s)), \nabla v(s)\Big)\,ds.  \ee
Применяя неравенство Коши-Буняковского, неравенство Коши и (5.4), мы получим:
\be\big| \int\limits_0^T\Big(f(s,x,v(s),\tau(s)), \nabla v(s)\Big)\,ds\big| \leq \|\widetilde f\|_{L_2((0,T)\times \Omega)}\|\nabla v\|_{L_2((0,T)\times \Omega)} $$ $$\leq  \frac 1 {2\Gamma_0}\|\widetilde f\|^2_{L_2((0,T)\times \Omega)}+  \frac{\Gamma_0} 2\|\nabla v\|^2_{L_2((0,T)\times \Omega)}.\ee

Аналогично,
\be\big| \int\limits_0^T\Big(g(s,x,v(s),\tau(s)), \nabla \tau(s)\Big)\,ds\big| \leq  \frac 1 {2\Gamma_0}\|\widetilde g\|^2_{L_2((0,T)\times \Omega)}+  \frac{\Gamma_0} 2\|\nabla \tau\|^2_{L_2((0,T)\times \Omega)}.\ee

Неравенства (5.11) -- (5.13) влекут
\be
 \Gamma^2\varepsilon\int\limits_0^T \|v(s)\|^2_2\,ds+ \varepsilon\int\limits_0^T\|\tau(s)\|^2_X\,ds + \frac {\lambda\Gamma_0\Gamma^2}2\int\limits_0^T \|\nabla v(s)\|^2\,ds+ \frac {\lambda\Gamma_0}2\int\limits_0^T\|\nabla\tau(s)\|^2\,ds $$ $$ \leq \frac {\Gamma^2} {2\Gamma_0}\|\widetilde f\|^2_{L_2((0,T)\times \Omega)}+\frac 1 {2\Gamma_0}\|\widetilde g\|^2_{L_2((0,T)\times \Omega)},  \ee
а потому и (5.6). $ \Box $

Дальнейший ход доказательства теоремы 5.1 (оценка производных, разреши\-мость вспомогательной задачи и предельный переход) проводится по аналогии с соот\-ветствующими рассуждениями из пункта 4.

\end {document}